\documentclass[a4paper,11pt,reqno,noindent]{amsart}
\usepackage[centertags]{amsmath}
\usepackage{mathtools}
\usepackage{amsfonts,amssymb,amsthm,dsfont,cases,amscd,esint,enumerate}
\usepackage[T1]{fontenc}
\usepackage[english]{babel}
\usepackage{psfrag}
\usepackage{newlfont,accents}
\usepackage{color}
\usepackage[body={15cm,21.5cm},centering]{geometry} 
\usepackage{fancyhdr}
\usepackage[font=small,labelfont=bf]{caption}

\usepackage{enumerate}

\mathchardef\emptyset="001F

\theoremstyle{plain}
\newtheorem{theorem}{Theorem}

\newtheorem{corollary}{Corollary}
\newtheorem{proposition}{Proposition}

\theoremstyle{definition}
\newtheorem{definition}{Definition}

\theoremstyle{remark}
\newtheorem{remark}{Remark}

\newcommand{\e}{\varepsilon}

\newcommand{\R}{\mathbb R}
\newcommand{\Z}{\mathbb Z}

\newcommand{\loc}{{\operatorname{loc}}}

\newcommand{\Id}{\operatorname{Id}}

\newcommand{\per}{{\operatorname{per}}}

\newcommand{\step}[1]{\noindent \textit{Step} #1.}

\newcommand{\RR}{{\mathbb{R}}}

\newcommand{\re}{\mathbb{R}}

\newcommand{\po}{\partial\Omega}

\usepackage{color}
\usepackage[colorlinks,citecolor=black,urlcolor=black]{hyperref}

\usepackage[normalem]{ulem}

\begin{document}
\title{Periodic homogenization and harmonic measures}

\begin{abstract}
Since the seminal work of Kenig and Pipher, the Dahlberg-Kenig-Pipher (DKP) condition on oscillations of the coefficient matrix became a standard threshold in the study of absolute continuity of the harmonic measure with respect to the Hausdorff measure on the boundary. It has been proved sufficient for absolute continuity in the domains with increasingly complex geometry, and known counterexamples show that in a certain sense it is necessary as well. 
\\
In the present note, we introduce into the subject ideas from homogenization theory to exhibit a new class of operators for which the elliptic measure is well-behaved, featuring the coefficients violating the DKP condition, and on the contrary, oscillating so quickly, that the homogenization takes place.
\end{abstract}

\author[G. David]{Guy David}
\address[Guy David]{Universit\'e Paris-Saclay, Laboratoire de Math\'{e}matiques d'Orsay, 91405, France}
\email{Guy.David@universite-paris-saclay.fr}
\author[A. Gloria]{Antoine Gloria}
\address[Antoine Gloria]{Sorbonne Universit\'e, Universit\'e Paris Cit\'e, CNRS, Laboratoire Jacques-Louis Lions, LJLL, F-75005 Paris, France  \& Universit\'e Libre de Bruxelles, D\'epartement de Math\'ematique, 1050~Brussels, Belgium}
\email{antoine.gloria@sorbonne-universite.fr}
\author[S. Qi]{Siguang Qi}
\address[Siguang Qi]{Sorbonne Universit\'e, Universit\'e Paris Cit\'e, CNRS, Laboratoire Jacques-Louis Lions, LJLL, F-75005 Paris, France}
\email{siguang.qi@sorbonne-universite.fr}
\author[S. Mayboroda]{Svitlana Mayboroda}
\address[Svitlana Mayboroda]{ETH Z\"urich, R\"amistrasse 101, 8092 Z\"urich Switzerland; School of Mathematics, University of Minnesota, 206 
Church St SE, Minneapolis, MN 55455 USA}
\email{svitlana.mayboroda@math.ethz.ch}

\date{}

\maketitle


\section{Introduction}

\subsection{History and motivation}

In this short article, we consider elliptic measures associated with operators of the form $L=-\nabla \cdot A \nabla$ in a domain $\Omega\subset \RR^{N+1}$. The last couple decades brought an incredible wealth of results pertaining to absolute continuity of the elliptic measure associated to $L$ under increasingly general conditions on the geometry of the domain. In this note, however, we will be content with $\Omega=\RR^{N+1}_+$ and instead concentrate on the coefficients of the equation.

In this direction, the two main threads of results pertain to some control of oscillations of the matrix $A$. Roughly speaking, either $A$ has to be independent of the transversal direction to the boundary \cite{JerisonKenig81}, or it should be quantitatively close to constant coefficient matrices. The latter condition, and the corresponding results in half-space, go back to Kenig and Pipher \cite{KP} where they referred to the underlying question as Dahlberg's conjecture, and at this point it has become a standard threshold for solvability of boundary value problems, referred, for short, as a DKP condition in most of the literature. The DKP condition is, in a certain sense, necessary and sufficient: one can find the corresponding counterexamples in \cite{CFK, FKP, Modica-Mortola, poggi}. 

Recent years have brought slightly weaker admissible variations of DKP \cite{DLM21}, however, one way or the other, these conditions amount to controlling quantitatively the (local) distance from $A$ to a constant matrix. There is however another way for an elliptic operator to be close to a constant-coefficient operator: in the resolvent sense rather than in the pointwise sense for $A$. This relates to the field of homogenization. The emerging operators have features opposite to DKP: they oscillate violently, so much so, that homogenization takes place. We do not strive for the biggest generality in the present note, and instead, present a case study, showing that the control of resolvents brought by the homogenization techniques is also sufficient for solvability of the boundary problems, absolute continuity of the harmonic measure, and related properties. 

We would like to point out that a different direction of employing homogenization in this subject has been tackled in the  series of papers \cite{MR2818717,MR2743875,MR2928140} by Kenig and Shen, and Lin. Homogenization is used  in \cite{MR2818717,MR2743875} to get scale-invariant estimates once small-scale estimates are established using DKP, whereas we use homogenization here as a proxy of DKP conditions to obtain the small scale estimates. In other words, Kenig, Lin, and Shen homogenize at large scales assuming that the small scales are under control for other reasons. Here, we homogenize at small scales. This requires delicate gluing arguments as our homogenization happens within individual Whitney cubes and the global correctors do not exist\footnote{A natural condition on $A$ would be fractal invariance towards the boundary, in which case one would get all scales at once. The difficulty is that there is an obstruction to constructing correctors with the same fractal invariance as $A$, which prevents the development of a global homogenization theory in this setting.}. 
%

Consider the  Dirichlet problem 
\begin{equation}\label{e.dir}
\begin{aligned}
-\nabla \cdot  A \nabla  u &=0 \quad  \textup{ in }\Omega,\\
u&=f \quad  \textup{ on }\partial \Omega,
\end{aligned}
\end{equation}
where $\Omega=\RR^{N+1}_+$ and $A$ is a  symmetric elliptic matrix with bounded measurable coefficients. For every $f\in C_0^\infty(\po)$ there exists a unique weak solution $u\in W^{1,2}(\Omega)$, continuous to the boundary. The solution can be represented by means of the corresponding harmonic measure.  In the half-space (and much bigger generality) absolute continuity of harmonic measure is equivalent to the solvability of the Dirichlet \eqref{e.dir} problem with singular, $L^p$, data, for some $p<\infty$, and equivalent to the Carleson measure estimates on the solutions below. It is convenient for us to focus on this latter reformulation, and the main theorem will be formulated in these terms.

For all $R>0$ and $x\in \R^N$, we define the Carleson tent of side-length $R$ centered at $x$
\[T_R(x)= \Bigl(x+(-\tfrac R2,\tfrac R2)^N\Bigr)  \times(0,R).\]  
By definition, tents are nested in the sense that $T_R(x) \subset T_{R'}(x)$ for all $R'\ge R$.
\begin{definition}
We say that a measure $\mu$ on $\R^N\times \R_+$ is a Carleson measure if there exists a constant $C>0$ such that for any $R>0$ and $x\in\re^N$,
\begin{equation}\label{def.carls}
\mu( T_R(x) )\le C R^N.
\end{equation}
\end{definition}
In \cite{KP}, Kenig and Pipher proved the absolute continuity of harmonic measure under the condition that $\Omega \ni Z \mapsto \tilde\alpha^2(Z) \frac{d Z}{\delta(Z)}$ is a Carleson measure, where 
$$\tilde\alpha(Z):= \delta(Z) \sup_{Y\in B_{\frac12 \delta(Z)}(Z)}|\nabla A(Y)|.$$
Here and throughout the paper $\delta(Z):=\textup{dist}(Z,\partial\Omega)$. One could also assume a slightly weaker condition below.
\begin{definition}\label{defDKP} We say that the matrix $A$ satisfies the DKP condition if 
$\alpha^2(Z) \frac{d Z}{\delta(Z)}$ is Carleson, for 
$$\alpha(Z):=  \sup_{Y,Y' \in B_{\frac12 \delta(Z)}(Z)} |A(Y)-A(Y')|.$$
\end{definition}
\begin{proposition}\label{Prop.DKP-Carls}\cite{KP, DPP} Assume that the matrix $A$ satisfies the DKP condition (cf. Definition~\ref{defDKP}). Then for any bounded solution $u$ to the boundary value problem \eqref{e.dir}, $t|\nabla u|^2\,dxdt$ is a Carleson measure: there exists a constant $C$ such that for any tent $T_R(x)$,
\begin{equation}\label{e.Carls}
\int_{T_R(x)} t|\nabla u(y,t)|^2 \,dydt\,\le\, C \|f\|_{L^\infty(\RR^N)}^2 R^N.
\end{equation}
\end{proposition}

This result is sharp, in the sense that 
\begin{theorem}\cite[Theorem~4.11]{FKP}\label{th:FKP}
Suppose that $\beta(Z) > 0$ is given on $\RR^2_+$, and satisfies the doubling condition: $\beta(Y) \le C\beta(Z)$ for all $Y \in B_{\frac12\delta(Z)}(Z)$. Assume that $\beta^2(Z) \frac{dZ}{\delta(Z)}$ is not a Carleson measure in  $[0,1]^2$.
Then there exists an elliptic operator $L = -\nabla \cdot A \nabla$ on $\RR^2_+$ such that, for $a(Z) :=
\sup_{Y,Y' \in B_{\frac12 \delta(Z)}(Z)} |A(Y) - A(Y') |$,
\begin{enumerate}
\item For any tent $T_R$, 
$$R^{-1}\int_{T_R}a^2(Z) \frac{dZ}{\delta(Z)}\le C\Big[R^{-1}\int_{T_{2R}}\beta^2(Z)\frac{dZ}{\delta(Z)} +1\Big].$$
\item The elliptic measure $\omega_L$ is not an $A_\infty$ measure on $[0, 1]$.
\end{enumerate}
\end{theorem}

In view of this and related counterexamples in \cite{CFK, Modica-Mortola, poggi} the theory rarely reaches beyond the DKP scenario. Here we present an alternative class of operators, failing DKP, for which the harmonic measure is well-behaved.

\subsection{Periodic homogenization}\label{sec:hom}

In Definition~\ref{defDKP}, the local oscillation of the coefficients monitors how close the coefficients are to a constant matrix in a pointwise sense. We recall here how closeness can be measured in the homogenization sense.

\medskip

Let $0<\lambda \le 1$ be a fixed ellipticity constant.
Let $A$ be a measurable 1-periodic in each direction symmetric matrix field that satisfies $\lambda \Id \le A \le \Id$.
Let $f \in L^2(\R^{N+1})^{N+1}$ be some forcing term and, for $0<\e \ll 1$, consider the solution $u_\e \in \dot H^1(\R^{N+1}):=\{v\in H^1_\loc(\R^{N+1})\,|\,\int_{\R^{N+1}} |\nabla v|^2 <\infty\}/\R$ of
\[
-\nabla \cdot A(\tfrac \cdot \e) \nabla u_\e =\nabla \cdot f.
\]
Homogenization aims at characterizing the large-scale behavior of $\nabla u_\e$ (see \cite[Chapter~1]{JKO-94}).
Indeed, at large scales (that is, scales much larger than $\e)$, $\nabla u_\e$ looks like the solution $\nabla \bar u$ of the equation 
\[
-\nabla \cdot \bar A \nabla \bar u =\nabla \cdot f,
\]
where $\lambda \Id \le \bar A\le \Id$ is a constant symmetric matrix depending only $A$ (and not on $f$).
For a gentle introduction to homogenization and the recent quantification thereof, we refer the reader to \cite{MR4648256}. In the rest of this article, we use the notation $\lesssim$ for $\le C\times$, where $C$ is universal constant (that depends on $N$, $\lambda$) which may vary from line to line.

\medskip

Before we go into detail, let us illustrate the type of convergence we are interested in on an elementary one-dimensional exemple, and consider the ODE on $(0,1)$ with $1$-periodic coefficient $a$
\[
(a(\tfrac \cdot \e) u'_\e)'(x)=f'(x), \quad u_\e(0)=u_\e(1)=0.
\]
The solution is explicitly given by
\begin{equation*}
u_\e(x)= \int_0^x \frac1{a(\frac y\e)}f(y)dy-\int_0^x \frac{1}{a(\frac y\e)}dy\bigg(\int_0^1 \frac1{a(\frac y\e)}dy\bigg)^{-1} \int_0^1 \frac1{a(\frac y\e)}f(y)dy,
\end{equation*}
which\footnote{Here we use that a rescaled periodic function converges weakly to its periodic average.} converges pointwise to 
\begin{equation*}
\bar u(x)= \int_0^x \tfrac1{\bar a}f(y)dy-x \int_0^1 \tfrac1{\bar a}f(y)dy, \quad \text{ where }\bar a = \Big(\int_0^1 \tfrac1{a(y)}dy\Big)^{-1},
\end{equation*}
the unique solution of 
\[
(\bar a \bar u')'(x)=f'(x), \quad \bar u(0)=\bar u(1)=0.
\]
Moreover, $\|u_\e-\bar u\|_{L^2(0,1)} \le C\e$. This entails that $(-\nabla \cdot a(\frac \cdot \e) \nabla)^{-1}$
is close to  $(-\nabla \cdot \bar a \nabla)^{-1}$ as $0<\e \ll 1$, which we shall use (albeit in a much more quantified form for gradients) as a replacement for the measure $\alpha$  in the DKP condition.

\medskip

In higher dimension, homogenization is not explicit. In order to define $\bar A$, we need to introduce the so-called correctors.
\begin{definition}\label{def:hom}
For all $1\le i\le N+1$, the corrector $\phi^i$ in direction $e_i$ (where $\{e_i\}_{1\le i\le N+1}$ denotes the canonical basis of $\R^{N+1}$) is 
the unique periodic solution in  $H^1_\per(Q)$ (the closure of smooth $Q$-periodic functions with average zero in $H^1(Q)$) of 
\[
-\nabla \cdot A(\nabla \phi^i +e_i)=0.
\]
The homogenized matrix $\bar A$ in direction $e_i$ is then given by
\[
\bar A e_i := \fint_Q A(\nabla \phi^i+e_i).
\]
We shall also need the flux corrector $\sigma^i$, defined in direction $e_i$ as the $(N+1) \times (N+1)$ skew-symmetric matrix field whose entries are the unique periodic solutions in $H^1_\per(Q)$ of
\[
-\Delta \sigma^{ijk} = \partial_j q^{ik}-\partial_k q^{ij},
\]
where $q^i:=A(\nabla \phi^i+e_i)-\bar A e_i$ and $q^{ik}=e_k \cdot q^i$.
By construction, we have $\nabla \cdot \sigma^i=q^i$, where the divergence is taken with respect to the last index, that is, $(\nabla \cdot \sigma^i)^j:=\sum_{k=1}^{N+1}\partial_k \sigma^{ijk}$.
\end{definition}
Homogenization is also quantitative at the level of gradients.
However, as one can already see in dimension $d=1$, $u_\e'$ has small-scale oscillations of order 1
and wavelength $\e$ whereas $\bar u'$ has no small-scale oscillations. To approximate $\nabla u_\e$ by $\nabla \bar u$, one needs to ``reconstruct'' these oscillations, which we do using  correctors $\{\phi^i\}_i$. To this aim, we introduce the two-scale expansion (we implicitly sum over $i$)
\begin{equation}\label{e.2s0} 
u_\e^{2s}:=\bar u + \e \phi^i (\tfrac \cdot \e) \partial_i \bar u.
\end{equation}
By using the flux corrector, we can show that the error $z_\e:=u_\e-u_\e^{2s}$ satisfies an equation (we implicitly sum over $i$)
\begin{equation}\label{e.2stim}
-\nabla \cdot A(\tfrac \cdot \e) \nabla z_\e = \e \nabla \cdot \Big((A \phi^i-\sigma^i)(\tfrac \cdot \e) \nabla \partial_i \bar u\Big)
\end{equation}
in $\RR^{n+1}$.
In particular, this entails the energy estimate
\[
\int_{\R^{N+1}} |\nabla z_\e|^2 \lesssim \e^2 \int_{\R^{N+1}} \sup_{i,Q} |(\phi^i,\sigma^i)|^2 |\nabla^2 \bar u|^2.
\]
The gain in the right-hand side comes from the regularity theory for the homogenized operator $-\nabla \cdot \bar A \nabla$ (which has constant coefficients): whereas $\nabla u_\e$ cannot be better than $L^2 \cap L^\infty(\R^{N+1})$, 
$\nabla \bar u$ is  in $H^1(\R^{N+1})$ provided $f \in L^2(\R^{N+1})$. This makes the right-hand side small, 
after dealing with the supremum by De Giorgi-Nash-Moser, which entails
\begin{equation}\label{e.6}
\int_Q \left(|\nabla \phi^i|^2 + |\nabla \sigma^i|^2\right)+\sup_Q( |\phi^i |+|\sigma^i|) \lesssim 1.
\end{equation}
Since \eqref{e.2stim} is at the core of our main result to come, we display the standard proof.
By scaling, we may assume that $\e=1$ (and we skip the subscript $\e$).
By definition of $z$ and adding and subtracting $\bar A \nabla \bar u$, 
\begin{equation}\label{e.defE}
A \nabla z =  A \nabla u-\bar A \nabla \bar u -\big(A(e_i+\nabla \phi^i)-\bar A e_i\big) \partial_i \bar u - A \phi^i \nabla \partial_i \bar u.
\end{equation}
Using that $\nabla \cdot \sigma^i=A(e_i+\nabla \phi^i)-\bar A e_i$
and $b \nabla \cdot c = \nabla \cdot (bc)-c \cdot \nabla b$, we may reformulate~\eqref{e.defE} as
\begin{equation}\label{e.flux0}
A   \nabla  z =  A   \nabla u - \bar A  \nabla \bar u 
-  (A \phi^i-\sigma^i)    \nabla  \partial_i \bar u
- \nabla \cdot (\sigma^i  \partial_i \bar u).
\end{equation}
We then recover \eqref{e.2stim} by using the following property: for all functions $w$,
\begin{equation}\label{e.divfree}
\nabla \cdot \big(\nabla \cdot (\sigma^i w )\big) \equiv 0.
\end{equation}
This is a simple consequence of the skew-symmetry of $\sigma^i$, as one can better see coordinate-wise. Indeed, $\nabla \cdot \big(\nabla \cdot (\sigma^i w )\big)=\partial_j (\partial_k (\sigma^{ijk} w))= (\partial^2_{jk} \sigma^{ijk})w+\sigma^{ijk} \partial_{jk}^2 w+\partial_k \sigma^{ijk}  \partial_j w+\partial_j \sigma^{ijk}  \partial_k w$. The first two terms obviously vanish by skew-symmetry, whereas the last two terms (which do not vanish individually) have opposite signs.

\medskip

In the upcoming subsection, rather than working on the whole space like above, we shall consider a localized version of \eqref{e.2s0}. Consider a smooth cut-off function $\chi \in C^\infty_c (\R^{N+1})$. The localized two-scale expansion then takes the form
\begin{equation}\label{e.2s} 
u_\e^{2s}:=\bar u + \e \phi^i (\tfrac \cdot \e) \chi \partial_i \bar u.
\end{equation}
A similar calculation as for \eqref{e.flux0} then yields (see for instance \cite[Proof of Proposition~1]{MR4103433} for details) for $z_\e:=u_\e-u_\e^{2s}$
\begin{multline}\label{e.flux}
A(\tfrac \cdot \e)  \nabla  z_\e =  A(\tfrac \cdot \e)  \nabla u_\e- \bar A  \nabla \bar u 
- \e (A \phi^i-\sigma^i) (\tfrac \cdot \e)  \nabla (\chi \partial_i \bar u)
-(1-\chi) (A(\tfrac \cdot \e)-\bar A) \nabla \bar u
\\
-\e \nabla \cdot (\sigma^i (\tfrac \cdot \e)\chi \partial_i \bar u),
\end{multline}
and
\begin{equation}\label{e.z}
-\nabla \cdot A(\tfrac \cdot \e) \nabla z_\e \,=\, \nabla \cdot \big( \e(A\phi^i-\sigma^i)(\tfrac \cdot \e) \nabla (\chi \partial_i \bar u)+(1-\chi) (A(\tfrac \cdot \e)-\bar A)\nabla \bar u)\big).
\end{equation}

\subsection{Main result: elliptic measures meet homogenization}\label{sec:main}

Consider the half space $\Omega=\R^{N+1}_+$.
Our aim is to show that homogenization can be used instead of DKP conditions when one gets close to the boundary $\partial \Omega:=\R^N \times \{0\}$. The main idea is to let the coefficients $A$ oscillate so fast when one gets close to $\partial \Omega$ that homogenization takes place fast enough to ensure that for any bounded solution $u$ to the boundary value problem \eqref{e.dir} $t|\nabla u|^2\,dxdt$ is a Carleson measure -- at least when restricted on small tents.

\medskip

We start by defining precisely the class of oscillating coefficients that we shall consider.
Since we are mostly interested in the local behavior close to $\partial \Omega$, we set for simplicity $A|_{t> 1}=A_\infty$, where $A_\infty$ is a constant matrix.
To describe the oscillations for $t<1$, we define Whitney boxes of the generation $k<0$ as the family of cubes $\{W_{kj}\}_{j\in \Z^N}$ of side-length $2^k$
\[
W_{kj}:=\Bigl(2^k j+\tfrac12[-2^{k},2^{k})^N\Bigr)\times [2^k,2^{k+1}).
\]
These boxes have the property that their side-length is comparable to their distance to the boundary $\partial \Omega$.
On each such cube we consider a $Q$-periodic matrix field $A_{kj}$ (as in Section~\ref{sec:hom})
and set $A|_{W_{kj}} := A_{kj}(\tfrac{\cdot}{2^k \e_{kj}})$,
where $0<\e_{kj}\le 1$ is a small parameter to be chosen later (and that will drive the smallness of the homogenization error).
All in all, we thus define
\begin{equation}\label{e.defA}
A := A_\infty \mathds{1}_{t\ge 1}+\sum_{k <0} \sum_{j \in \Z^N} A_{kj}(\tfrac{\cdot}{2^k \e_{kj}}) \mathds1_{W_{kj}}.
\end{equation}
The field $A|_{t < 1}$ is locally periodic (on each Whitney cube). Its period goes to zero as one approaches $\partial \Omega$, even faster than the size of the Whitney cubes decreases. More precisely, in the Whitney box $W_{kj}$, $A$ is $2^k \e_{kj}$-periodic -- and $A$ has $\e_{kj}^{-1}$ periods per dimension in $W_{kj}$.

To each Whitney box $W_{kj}$, we can associate the (constant) homogenized matrix $\bar A_{kj}$ of the matrix field $A_{kj}$ via Definition~\ref{def:hom}, and we set 
\begin{equation}\label{e.barA}
\bar A := A_\infty \mathds{1}_{t\ge 1}+\sum_{k<0} \sum_{j \in \Z^N} \bar A_{kj}  \mathds1_{W_{kj}},
\end{equation}
which is a piece-wise constant matrix field on $\R^N \times \R_+$.
For the sake of illustration, we display a cartoon of $A$ on Figure~\ref{fig1}.

\medskip

\begin{center}
\includegraphics[width=0.70\linewidth]{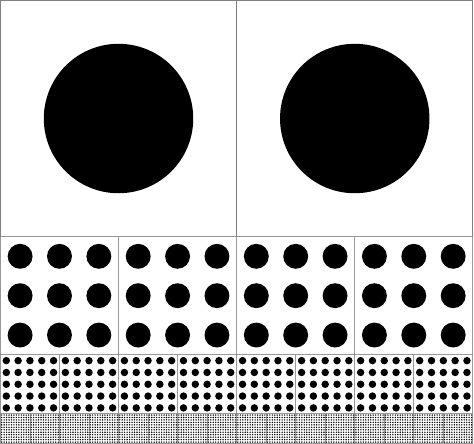}

\vspace{-10pt}

\captionof{figure}{Sketch of $A$ for $t\in [\frac1{16},1]$
 (that is, 4 generations) in case of locally periodic inclusions: the periodic structure gets finer as one gets closer to the boundary. } \label{fig1}

\vspace{10pt}
\end{center}

\medskip

Our main result is as follows.
\begin{theorem}\label{theor1} Let $\Omega=\RR^{N+1}_+.$
For all $0< \lambda \le 1$, denote by $2p>2$ the Meyers' exponent\footnote{which depends on $N,\lambda$ and satisfies $\lim_{\lambda \to 1} p(N,\lambda)=\infty$} discussed near \eqref{Meyers}. 
Let $A$ and $\bar A$ be given by \eqref{e.defA} and \eqref{e.barA}, respectively. Assume that there exists a constant $C>0$ such that for all $f\in H^{1/2}(\partial \Omega)\cap L^\infty (\partial\Omega),$ the unique weak solution $\bar u$ to 
\begin{equation}\label{e.hom}
\begin{aligned}
-\nabla \cdot  \bar A \nabla  \bar u &=0 \quad  \textup{ in }\Omega,\\
\bar u&=f \quad  \textup{ on }\partial \Omega,
\end{aligned}
\end{equation}
satisfies
\begin{equation}\label{eqCM1}
\sup_{x\in \RR^N, \, R>0} R^{-N}\int_{T_R(x)} t|\nabla \bar u|^2 \leq C \|f\|_{L^\infty(\RR^N)}. 
\end{equation}
If for $k<0$ and $j \in \Z^N$ we have $\e_{kj}\lesssim 2^{\alpha(p) k} $ with $\alpha(p):=\frac{3p-1}{2(p-1)}$,
then the weak solution $u$ to the boundary value problem \eqref{e.dir} also satisfies
$$\sup_{x\in \RR^N, \, R>0} R^{-N}\int_{T_R(x)} t|\nabla  u|^2 \leq C \|f\|_{L^\infty(\RR^N)}, $$
with a constant $C$ that depends on the constant in \eqref{eqCM1} and our choice of $\e_{kj}$ but not on $f$.

%
%
\end{theorem}
\begin{remark}
Comments are in order:
\begin{itemize}
\item The Carleson measure estimates on the solutions assumed for $-\nabla \cdot \bar A \nabla$ hold by Proposition~\ref{Prop.DKP-Carls} if we assume that $\bar A$ (which is piece-wise constant) satisfies the DKP condition of Definition~\ref{defDKP}.
\item The proof is local in the sense that homogenization controls the solution on $T_R \cap \{t<2\}$ -- the assumption that $A|_{t>1}$ be constant is only used to get the Carleson estimates on tents of size $R\ge 2$.
\item One could in principle mix the present approach and the approach of \cite{MR2743875} by periodizing $A|_{0<t<1}$ in the vertical direction $t$, and therefore use homogenization both at small and large scales.
\end{itemize}
\end{remark}
Let us give a concrete example, and consider the case when $A_{kj} = A_\per$ is a non-constant periodic field that does not depend on $k,j$ (yet, the period of $A|_{W_{kj}}$ strongly depends on $k,j$) and choose $A_\infty=\bar A_\per$ (the homogenized coefficient of $A_\per$). Then $\bar A \equiv \bar A_\per$ is constant, so that  $-\nabla \cdot \bar A \nabla$ has constant coefficients and Theorem~\ref{theor1} applies. Yet, $A$ trivially violates the DKP condition.
Indeed, for all $Z\in \R^N \times (0,\frac12)$,
\[
\alpha(Z)=  \sup_{Y,Y' \in B_{\frac12 \delta(Z)}(Z)} |A(Y)-A(Y')|\sim 1,
\]
and $\int_{T_{\frac 12}(0)}\alpha(Z)^2\frac{dZ}{\delta(Z)}=\infty$. 
This example can be seen as a positive counterpart (due to homogenization) to the counter-example of Theorem~\ref{th:FKP}.

\medskip 

We conclude with the observation that in the case when the homogenized matrix field $\bar A$ has better properties
than assumed for Theorem~\ref{theor1}, one can control the homogenization error slightly better.
\begin{corollary}\label{cor1}
If $\bar A$ is a laminate in the sense that $x \mapsto \bar A(x)$ only depends on $x_{N+1}$, then 
Theorem~\ref{theor1} holds under the weaker condition $\e_{kj}\lesssim 2^{\frac32 k} $. If $\bar A$ is constant, then Theorem~\ref{theor1} holds under the even weaker condition $\e_{kj}\lesssim 2^{ k} $.
\end{corollary}

\section{Proof of Theorem~\ref{theor1}}

Recall that $f\in H^{\frac12}(\partial \Omega) \cap L^\infty(\partial \Omega)$.
Set $S_1:=\{ (x,t) \in \Omega \,|\, t < 1\}$ and $S_2:=\{ (x,t) \in \Omega \,|\, t < 2\}$.
The strategy of the proof is as follows.
By our simplifying assumption that $A|_{t\ge 1}=A_\infty$ is constant, 
we first observe that  for $R\ge 2$,
\begin{equation}\label{e.up}
\int_{T_R\backslash S_1}(t-1)|\nabla u|^2\lesssim R^N \|f\|_{L^\infty}^2.
\end{equation}
Indeed, since $A$ is constant in $\R^N \times (1,+\infty)$, Proposition~\ref{Prop.DKP-Carls} applied to $u$ on $\R^N \times (1,+\infty)$ yields for all $R\ge 2$
\[
\int_{T_R\backslash S_1}(t-1)|\nabla u|^2\lesssim R^N\|u\|_{L^\infty(\R^N\times \{t=1\})}^2.
\]
Therefore, \eqref{e.up} follows from the maximum principle in form of $\|u\|_{L^\infty(\R^N\times \{t=1\})}\le \|f\|_{L^\infty}$.

\medskip

Hence it remains to prove that for all $R>0$, 
\begin{equation} \label{e.down}
\int_{T_R \cap S_2} t|\nabla u|^2 \lesssim R^{N}\|f\|_{L^\infty}^2.
\end{equation}
The proof of \eqref{e.down} takes advantage of homogenization and our assumption on $\bar A$.
Recall that $\bar u$ is the solution of \eqref{e.hom}.
We shall approximate $u$ by local two-scale expansions based on $\bar u$, and need further notation:
\begin{itemize}
\item $\hat W_{kj}$ denotes the union of the Whitney cubes that touch $W_{kj}$.
\item $\dot W_{kj}^\eta$ denotes the cube of sidelength $(1-\eta) 2^{k}$ centered at the center of $W_{kj}$, where $0<\eta \le \frac12$ measures the width of the boundary layer in the Whitney cube in a scale-invariant way.
We write $\dot W_{kj}:=\dot W_{kj}^{\eta_{kj}}$ where $\eta_{kj}$ will be chosen in the proof. 
\item $\chi_{kj}\in C_c^\infty(W_{kj})$ is a smooth non-negative cut-off function for $\dot W_{kj}$ in $\dot W_{kj}^{\eta_{kj}/2}$, that is, $\chi_{kj}\equiv1$ in $\dot W_{kj}$, $\chi_{kj}\equiv0$ in $W_{kj} \setminus\dot  W_{kj}^{\eta_{kj}/2}$, and $\sup (\chi_{kj} + 2^k \eta_{kj}|\nabla \chi_{kj}|)\lesssim 1$.
\end{itemize}
We now fix a tent $T_R(x)$ with $R>0$. Since the location of the tent $x$ plays no special rule in the proof, we will drop it in the following. Without loss of generality, we assume that $\|f\|_{L^\infty}=1$.
We define the local two-scale expansion
\begin{equation}\label{e.a02}
u^{2s}:=\bar u+ \sum_{k<0}\sum_{j\in \Z^N}  2^{k} \e_{kj}\chi_{kj} \phi_{kj}^i(\tfrac \cdot {2^{k}\e_{kj}}) \partial_i \bar u,
\end{equation}
where the sum on $k,j$ is taken over those pairs such that $W_{kj} \subset S_1$.  
Although it is not clear a priori, we have $u^{2s} \in H^1(\Omega)$.
Since  $f\in H^\frac12(\partial \Omega)$, $\bar u \in H^1(\Omega)$.
Hence, $u^{2s} \in L^2(\Omega)$ using the fact that  $|\chi_{kj}^i|\lesssim 1$ and \eqref{e.6}.
The fact that $\nabla u^{2s} \in L^2(\Omega)$ is a by-product of Step~1 below (by summing \eqref{e.H1!}
over $k<0$ and $j\in \Z^N$, using that $\nabla u^{2s}=\nabla \bar u$ for $t>1$, and $\bar u \in H^1(\Omega)$).  
As a consequence, note also that $u^{2s}|_{\partial \Omega}=f$ in $H^\frac12(\partial \Omega)$
(each term of the sum in \eqref{e.a02} is compactly supported in $\Omega$ and thus vanishes on $\partial \Omega$).
We shall obtain \eqref{e.down} as the combination of the following two estimates: on the one hand 
\begin{equation} \label{e.1}
\int_{T_R } t|\nabla u^{2s}|^2 \lesssim R^{N},
\end{equation}
and on the other hand  
\begin{equation} \label{e.2}
\int_{T_R \cap S_2} t |\nabla (u^{2s}-u)|^2 \lesssim R^{N}.
\end{equation}
We prove \eqref{e.1} and \eqref{e.2} in the upcoming two steps.

\medskip

\step1 Proof of~\eqref{e.1} provided $\e_{kj}\lesssim \eta_{kj}$. \\
By assumption, we have 
\begin{equation}\label{cald-3}
\int_{T_R} t|\nabla \bar u|^2\,\lesssim \,  R^N,
\end{equation}
to which we shall reduce the claim.
To this aim, we start with differentiating the formula~\eqref{e.a02} for $u^{2s}$, which yields on each $W_{kj}$ with $k<0$
\begin{equation}\label{u2s}
\nabla u^{2s} = \Big(e_i+ \chi_{kj} \nabla \phi_{kj}^i(\tfrac \cdot {2^{k} \e_{kj}})+2^{k}\e_{kj}  \phi_{kj}^i(\tfrac \cdot {2^{k}\e_{kj}}) \nabla\chi_{kj} \Big) \partial_i \bar u
+2^{k} \e_{kj} \chi_{kj} \phi_{kj}^i(\tfrac \cdot {2^{k}\e_{kj}}) \nabla \partial_i \bar u,
\end{equation}
whereas $\nabla u^{2s}=\nabla \bar u$ for $t \ge 1$. 

We start with the control of the first right-hand side term of \eqref{u2s}.
Since $\phi_{kj}^i$ is bounded (cf.~\eqref{e.6}) and provided $\e_{kj} \lesssim \eta_{kj}$
so that $|2^{k}\e_{kj} \nabla\chi_{kj}|\lesssim 1$, we have 
\[
\Big|(e_i+2^{k}\e_{kj} \nabla\chi_{kj} \phi_{kj}^i(\tfrac \cdot {2^{k}\e_{kj}}))\partial_i \bar u\Big|\, \lesssim\,|\nabla \bar u|,
\]
and it remains to treat the term involving $\nabla \phi_{kj}^i$, which is only controlled when locally squared-averaged (cf.~\eqref{e.6}). To this aim we shall use that $\bar u$ is $\bar A_{kj}$-harmonic in $W_{kj}$, to the effect that for all cubes $Q_{2^{k+1}\e_{kj}}$ of side-length $2^{k+1}\e_{kj}$ contained in $W_{kj}$ we have for the cube $Q_{2^{k}\e_{kj}}$ with the same center and half the side-length 
\[
\sup_{Q_{2^{k}\e_{kj}}} |\nabla \bar u|^2 \,\lesssim\, \fint_{Q_{2^{k+1}\e_{kj}}}|\nabla \bar u|^2.
\]
In particular, by~\eqref{e.6}, this entails
\[
\int_{Q_{2^{k}\e_{kj}}} |\nabla \phi_{kj}^i(\tfrac \cdot {2^{k} \e_{kj}})|^2 |\nabla \bar u|^2 \,\lesssim\,\int_{Q_{2^{k+1}\e_{kj}}}|\nabla \bar u|^2.
\]
Assuming that $\eta_{kj}\gtrsim \e_{kj}$,  the support of $\chi_{kj}$ can be covered by such cubes $Q_{2^{k}\e_{kj}}$, so that 
\begin{equation}\label{e.a01}
\int_{W_{kj}} \Big|\Big(e_i+ \chi_{kj} \nabla \phi_{kj}^i(\tfrac \cdot {2^{k} \e_{kj}})+2^{k}\e_{kj}  \phi_{kj}^i(\tfrac \cdot {2^{k}\e_{kj}}) \nabla\chi_{kj} \Big) \partial_i \bar u\Big|^2 \, \lesssim\, \int_{W_{kj}} |\nabla \bar u|^2.
\end{equation}

We now turn to the second right-hand term of  \eqref{u2s}, and need to control the second derivatives by the gradient,
which we do again by elliptic regularity using the fact that  $\bar A_{kj}$ is constant in $W_{kj}$.
By Caccioppoli's inequality applied to $\nabla \bar u$  with the cut-off function $\chi_{kj}^2$
\begin{equation}\label{cald-5}
\int_{W_{kj}} \chi_{kj}^2 |\nabla^2 \bar u|^2 \,\lesssim\, \big(\tfrac{1}{2^k\eta_{kj}}\big)^2\int_{W_{kj}} |\nabla \bar u|^2.
\end{equation}
By  \eqref{e.6}, this implies 
\[
\int_{W_{kj}}  \Big| 2^k\e_{kj} \chi_{kj} \phi_{kj}^i(\tfrac \cdot {2^{k}\e_{kj}}) \nabla \partial_i \bar u\Big|^2  
\,\lesssim\, \big(\tfrac{\e_{kj}}{\eta_{kj}}\big)^2\int_{W_{kj}} |\nabla \bar u|^2.
\]
In combination with \eqref{e.a01} this yields for all $k<0$, under the condition $\e_{kj}\lesssim \eta_{kj}$ which we retain,
\begin{equation}\label{e.H1!}
\int_{W_{kj}} |\nabla u^{2s} |^2   \,\lesssim\, \int_{W_{kj}} |\nabla \bar u|^2.
\end{equation}
Since $2^{k} \le t \le 2^{k+1}$ on $W_{kj}$, we can add the weight $t$ in each Whitney cube.
Using that $\nabla u^{2s}=\nabla \bar u$ for $t\ge 1$, the claim~\eqref{e.1} then follows from~\eqref{cald-3}.

\medskip

\step2 Proof of~\eqref{e.2}.\\
This is where homogenization comes into play, which allows us to establish~\eqref{e.2} in the (stronger\footnote{it is stronger because we upper bound the weight in~\eqref{e.2}, and then replace $T_R \cap S_2$ by $T_R$ -- yet we do not control the integral of $t|\nabla (u^{2s}-u)|^2$ on $T_R$.}) form of 
\begin{equation} \label{e.2b}
\int_{T_R} (R \wedge 1) |\nabla (u^{2s}-u)|^2 \lesssim R^{N}.
\end{equation}
On each Whitney cube $W_{kj}$ with $k<0$, the homogenization error $z:=u-u^{2s}\in H^1_0(\Omega)$ of the two-scale expansion satisfies, according to~\eqref{e.z}, the equation
\[
-\nabla \cdot A \nabla z\,=\,  \nabla \cdot f_{kj},
\]
with
\[
f_{kj}:=2^{k}\e_{kj}(A\phi_{kj}^i(\tfrac \cdot {2^{k} \e_{kj}})-\sigma^i_{kj} (\tfrac \cdot {2^{k} \e_{kj}})\nabla(\chi_{kj} \partial_i \bar u))
+ \mathds{1}_{W_{kj}}(1-\chi_{kj})(A-\bar A_{kj}) \nabla \bar u
\]
for $k<0$. This equation also holds for $k\ge 0$ with $f_{kj}\equiv 0$. Summing over $k \in \Z,j\in \Z^N$, we then get
\begin{equation}\label{e.4}
\begin{array}{rcll}
-\nabla \cdot A \nabla z&=&  \nabla \cdot f  & \textup{ in }\Omega,\\
z&=&0  & \textup{ on }\partial\Omega,
\end{array}
\end{equation}
with $f:=\sum_{k<0}\sum_{j\in \Z^N} f_{kj}$. 

With \eqref{e.4} at hand, using that $z$ vanishes on $\partial \Omega$, we can appeal to Caccioppoli's inequality with a cut-off function for $T_R$ in $T_{2R}$, to the effect that
\begin{equation}\label{e.fin0}
\int_{T_R} |\nabla z|^2 \,\lesssim \, \int_{T_{2R}} |f|^2 +R^{-2} \int_{T_{2R}} z^2.
\end{equation}
We estimate both terms on the right-hand side of  \eqref{e.fin0} separately, and start with the second one.

By the maximum principle $\sup |u_\e|, \sup |\bar u|\le 1$, and hence
\[
R^{-2} \int_{T_{2R}} \bar u^2+u_\e^2  \lesssim R^{N-1}.
\]
It remains to treat the term in $z$ involving the corrector.
By \eqref{e.6}, and using that $2^{k} \e_{kj} \lesssim R\wedge 1$ on $T_{2R}$,
\[
\int_{T_{2R}}\Big|\sum_{k<0} \sum_{j\in \Z^N} 2^{k} \e_{kj}\chi_{kj} \phi_{kj}^i(\tfrac \cdot {2^{k}\e_{kj}}) \partial_i \bar u \Big|^2\,\lesssim\, ( R\wedge 1)^2\int_{T_{2R} \cap S_1} |\nabla \bar u|^2.
\]
By Caccioppoli's inequality for $\bar u$ with a cut-off function for $T_{2R} \cap S_1$ in $T_{4R} \cap S_2$, followed by the maximum principle in form of $\sup |\bar u|\le 1$,
this yields
\[
\int_{T_{2R}}\Big|\sum_{k<0} \sum_{j\in \Z^N} 2^{k} \e_{kj}\chi_{kj} \phi_{kj}^i(\tfrac \cdot {2^{k}\e_{kj}}) \partial_i \bar u \Big|^2\,\lesssim\, \int_{T_{4R} \cap S_2} \bar u^2 \,\lesssim\, R^{N+1}.
\]
We have thus proved
\begin{equation}\label{e.fin1}
R^{-2} \int_{T_{2R}} z^2  \, \lesssim\,  R^{N-1}.
\end{equation}

We now control the first right-hand side term of \eqref{e.fin0}. 
A direct calculation yields 
\begin{eqnarray}
\int_{T_{2R}} |f|^2 &\lesssim& \sum_{kj} \int_{W_{kj}} (2^{k}\e_{kj})^2(|\nabla \chi_{kj}|^2 |\nabla \bar u|^2
+\chi_{kj}^2|\nabla^2 \bar u|^2)+(1-\chi_{kj})^2 |\nabla \bar u|^2\nonumber
\\
&\lesssim&  \sum_{kj} \int_{W_{kj}} (2^{k}\e_{kj})^2\chi_{kj}^2|\nabla^2 \bar u|^2
+\sum_{kj} \Bigl(1+(\frac{\e_{kj}}{\eta_{kj}})^2\Bigr)  \int_{W_{kj}\setminus \dot W_{kj}} |\nabla \bar u|^2,
\label{e.canbeimproved}
\end{eqnarray}
where the sum over $k,j$ runs over indices $k<0$ and $j\in \Z^N$  such that $W_{kj}\cap  T_{2R}\cap S_1 \ne \varnothing$. 
The first term on the right-hand side has been already controlled in \eqref{cald-5}. The second one is a boundary layer that has to be treated differently. Although $\bar A$ is locally constant, we do not (necessarily) have Lipschitz regularity at corners (cf.~the counterexample of \cite{MR1770930}), and we only have Meyers' estimate at our disposal. By H\"older's inequality with exponents $(p, \frac p{p-1})$ and Meyer's inequality (or reverse H\"older inequality for gradients, see e.g.~\cite{GM} and \cite{U}),
\begin{equation}\label{Meyers}
\int_{W_{kj}\setminus \dot W_{kj}} |\nabla \bar u|^2\, \lesssim \, (\eta_{kj}2^{k(N+1)})^\frac{p-1}{p} 
\Big(\int_{W_{kj}} |\nabla \bar u|^{2p}\Big)^\frac1p \, \lesssim\, \eta_{kj}^{\frac{p-1}{p}}\int_{\hat W_{kj}} |\nabla \bar u|^{2}.
\end{equation}
All in all, we thus have
\[
\int_{T_{2R}} |f|^2\,\lesssim\, \sum_{kj} \Bigl( (\frac{\e_{kj}}{\eta_{kj}})^2+\eta_{kj}^\frac{p-1}p\Bigr) \int_{\hat W_{kj}} |\nabla \bar u|^2.
\]
Given $\e_{kj}$ we optimize in $\eta_{kj}$ by choosing $\eta_{kj}=\e_{kj}^{\frac{2p}{3p-1}}$ (for which the proviso $\eta_{kj} \gtrsim \e_{kj}$ of Step~1 holds true), and obtain
\[
\int_{T_{2R} } |f|^2\,\lesssim\, \sum_{kj} \e_{kj}^{\frac{2(p-1)}{3p-1}} \int_{\hat W_{kj}} |\nabla \bar u|^2.
\]
By our choice of $\e_{kj}$, $\e_{kj}^{\frac{2(p-1)}{3p-1}}\lesssim t$ on $\hat W_{kj}$, and the assumption that $t|\nabla \bar u|^2$ is a Carleson measure entails
\begin{equation}\label{e.fin2}
\int_{T_{2R}} |f|^2\,\lesssim\, \int_{T_{4R} \cap S_2} t |\nabla \bar u|^2 \, \le \, \int_{T_{4R}  } t |\nabla \bar u|^2 \, \lesssim\,R^N.
\end{equation}
The desired estimate \eqref{e.2b} now follows from \eqref{e.fin0}, \eqref{e.fin1}, and \eqref{e.fin2}
in form of
\[
(R \wedge 1) \int_{T_R} |\nabla z|^2 \lesssim (R\wedge 1)R^N+(R\wedge 1)R^{N-1} \lesssim R^N.
\]

\section{Proof of Corollary~\ref{cor1}}
If $\bar A$ is a laminate (that is, only a function of $t$), then we have uniform Lipschitz regularity by \cite[Theorem~2]{MR853976}, to the effect that 
\begin{equation}\label{Lip}
\sup_{W_{kj}} |\nabla \bar u|^2 \,\lesssim\,  \fint_{\hat W_{kj}} |\nabla \bar u|^2,
\end{equation}
which is an upgrade of \eqref{Meyers} (which amounts to $p=\infty$, and therefore $\frac{2p}{3p-1}=\frac23$ in the optimization process). Notice that the condition that $\bar A$ is laminate is not necessary, we only need the Lipschitz estimate \eqref{Lip} which could be valid for other reasons.

\medskip

If $\bar A$ is a constant matrix, then $\bar u$ is harmonic, and we may also upgrade \eqref{cald-5}  to
\begin{equation}\label{cald-5+}
\sup_{W_{kj}} (|\nabla \bar u|^2+2^{2k}|\nabla^2 \bar u|^2) \,\lesssim\,   \fint_{\hat W_{kj}} |\nabla \bar u|^2.
\end{equation}
With \eqref{cald-5+} at hand, one may improve \eqref{e.canbeimproved} to
\begin{eqnarray*}
\int_{T_{2R}} |f|^2&\lesssim &    \sum_{kj} \int_{W_{kj}} (2^{k}\e_{kj})^2\chi_{kj}^2|\nabla^2 \bar u|^2
+\sum_{kj} \Bigl(1+(\frac{\e_{kj}}{\eta_{kj}})^2\Bigr)  \int_{W_{kj}\setminus \dot W_{kj}} |\nabla \bar u|^2 
\\
&\lesssim& \sum_{kj} \Bigl(\e_{kj}^2+\eta_{kj}+\frac{\e_{kj}^2}{\eta_{kj}}\Bigr)\int_{\hat W_{kj}} |\nabla \bar u|^2
\,\lesssim\, \sum_{kj}\e_{kj}\int_{\hat W_{kj}} |\nabla \bar u|^2
\end{eqnarray*}
with the choice $\eta_{kj} \sim \e_{kj}$. For $\e_{kj}\lesssim 2^k$, this entails \eqref{e.fin2} as well.

\section*{Acknowledgements}

AG acknowledges financial support from the European Research Council (ERC) under the European Union's Horizon 2020 research and innovation programme (Grant Agreement n$^\circ$~864066). GD is supported by Simons Foundation grant 601941, GD. SM is supported in part by the NSF grant DMS 1839077 and Simons Foundation grant 563916, SM.

\bibliographystyle{abbrv}
\bibliography{HomHarm-4}

\begin{thebibliography}{10}

\bibitem{CFK}
L.~A. Caffarelli, E.~B. Fabes, and C.~E. Kenig.
\newblock Completely singular elliptic-harmonic measures.
\newblock {\em Indiana Univ. Math. J.}, 30(6):917--924, 1981.

\bibitem{MR853976}
M.~Chipot, D.~Kinderlehrer, and G.~Vergara-Caffarelli.
\newblock Smoothness of linear laminates.
\newblock {\em Arch. Rational Mech. Anal.}, 96(1):81--96, 1986.

\bibitem{DLM21}
G.~David, L.~Li, and S.~Mayboroda.
\newblock Carleson measure estimates for the {G}reen function.
\newblock {\em Arch. Ration. Mech. Anal.}, 243(3):1525--1563, 2022.

\bibitem{DPP}
M.~Dindos, S.~Petermichl, and J.~Pipher.
\newblock The {$L^p$} {D}irichlet problem for second order elliptic operators
  and a {$p$}-adapted square function.
\newblock {\em J. Funct. Anal.}, 249(2):372--392, 2007.

\bibitem{FKP}
R.~A. Fefferman, C.~E. Kenig, and J.~Pipher.
\newblock The theory of weights and the {D}irichlet problem for elliptic
  equations.
\newblock {\em Ann. of Math. (2)}, 134(1):65--124, 1991.

\bibitem{GM}
M.~Giaquinta and L.~Martinazzi.
\newblock {\em An introduction to the regularity theory for elliptic systems,
  harmonic maps and minimal graphs}, volume~11 of {\em Appunti. Scuola Normale
  Superiore di Pisa (Nuova Serie) [Lecture Notes. Scuola Normale Superiore di
  Pisa (New Series)]}.
\newblock Edizioni della Normale, Pisa, second edition, 2012.

\bibitem{MR4648256}
A.~Gloria.
\newblock Homog\'en\'eisation stochastique: du qualitatif au quantitatif.
\newblock {\em Gaz. Math.}, (171):41--51, 2022.

\bibitem{MR4103433}
A.~Gloria, S.~Neukamm, and F.~Otto.
\newblock A regularity theory for random elliptic operators.
\newblock {\em Milan J. Math.}, 88(1):99--170, 2020.

\bibitem{JerisonKenig81}
D.~S. Jerison and C.~E. Kenig.
\newblock The {D}irichlet problem in nonsmooth domains.
\newblock {\em Ann. of Math. (2)}, 113(2):367--382, 1981.

\bibitem{JKO-94}
V.~V. Jikov, S.~M. Kozlov, and O.~A. Ole{\u \i}nik.
\newblock {\em Homogenization of differential operators and integral
  functionals}.
\newblock Springer-Verlag, Berlin, 1994.

\bibitem{MR1770930}
C.~Kenig, H.~Koch, J.~Pipher, and T.~Toro.
\newblock A new approach to absolute continuity of elliptic measure, with
  applications to non-symmetric equations.
\newblock {\em Adv. Math.}, 153(2):231--298, 2000.

\bibitem{MR2928140}
C.~E. Kenig, F.~Lin, and Z.~Shen.
\newblock Convergence rates in {$L^2$} for elliptic homogenization problems.
\newblock {\em Arch. Ration. Mech. Anal.}, 203(3):1009--1036, 2012.

\bibitem{KP}
C.~E. Kenig and J.~Pipher.
\newblock The {D}irichlet problem for elliptic equations with drift terms.
\newblock {\em Publ. Mat.}, 45(1):199--217, 2001.

\bibitem{MR2818717}
C.~E. Kenig and Z.~Shen.
\newblock Homogenization of elliptic boundary value problems in {L}ipschitz
  domains.
\newblock {\em Math. Ann.}, 350(4):867--917, 2011.

\bibitem{MR2743875}
C.~E. Kenig and Z.~Shen.
\newblock Layer potential methods for elliptic homogenization problems.
\newblock {\em Comm. Pure Appl. Math.}, 64(1):1--44, 2011.

\bibitem{Modica-Mortola}
L.~Modica and S.~Mortola.
\newblock Construction of a singular elliptic-harmonic measure.
\newblock {\em Manuscripta Math.}, 33(1):81--98, 1980/81.

\bibitem{poggi}
B.~Poggi.
\newblock Failure to slide: a brief note on the interplay between the
  kenig-pipher condition and the absolute continuity of elliptic measures.
\newblock 2019.
\newblock arXiv:1912.10115.

\bibitem{U}
K.~Uhlenbeck.
\newblock A new proof of a regularity theorem for elliptic systems.
\newblock {\em Proc. Amer. Math. Soc.}, 37:315--316, 1973.

\end{thebibliography}

\end{document}